\documentclass[12pt]{article}%
\usepackage{amsmath,amssymb,fullpage}
\usepackage{showlabels}
\usepackage[applemac]{inputenc}
\usepackage{amsmath}
\usepackage{amsfonts}
\usepackage{amssymb}
\usepackage{graphicx}%
\providecommand{\U}[1]{\protect\rule{.1in}{.1in}}

\newtheorem{example}{Example}[section]
\newtheorem{theorem}[example]{Theorem}
\newtheorem{definition}[example]{Definition}
\newtheorem{lemma}[example]{Lemma}
\newtheorem{problem}[example]{Problem}
\newtheorem{remark}[example]{Remark}

\numberwithin{equation}{section}

\def\1B{\text{1\!\!I}}

\begin{document}

\title{Stochastic Control of Memory Mean-Field Processes}
\author{Nacira A{\small GRAM}$^{1,2}$ and Bernt Ø{\small KSENDAL}$^{1,2}$}
\date{25 October 2017 \vskip0.3cm \emph{Dedicated to the memory of Salah-Eldin Mohammed}}
\maketitle

\footnotetext[1]{Department of Mathematics, University of Oslo, P.O. Box 1053
Blindern, N--0316 Oslo, Norway.\newline Email: \texttt{naciraa@math.uio.no,
oksendal@math.uio.no}}

\footnotetext[2]{This research was carried out with support of the Norwegian
Research Council, within the research project Challenges in Stochastic
Control, Information and Applications (STOCONINF), project number 250768/F20.}


\paragraph{Abstract}

By a memory mean-field process we mean the solution $X(\cdot)$ of a stochastic
mean-field equation involving not just the current state $X(t)$ and its law
$\mathcal{L}(X(t))$ at time $t$, but also the state values $X(s)$ and its law
$\mathcal{L}(X(s))$ at some previous times $s<t$. Our purpose is to study
stochastic control problems of memory mean-field processes.

\begin{itemize}
\item We consider the space $\mathcal{M}$ of measures on $\mathbb{R}$ with the
norm $|| \cdot||_{\mathcal{M}}$ introduced by Agram and Øksendal in
\cite{AO1}, and prove the existence and uniqueness of solutions of memory
mean-field stochastic functional differential equations.

\item We prove two stochastic maximum principles, one sufficient (a
verification theorem) and one necessary, both under partial information. The
corresponding equations for the adjoint variables are a pair of \emph{(time-)
advanced backward stochastic differential equations}, one of them with values
in the space of bounded linear functionals on path segment spaces.

\item As an application of our methods, we solve a memory mean-variance
problem as well as a linear-quadratic problem of a memory process.
\end{itemize}



\section{Introduction}

In this work we are studying a general class of controlled memory mean-field
stochastic functional differential equations (mf-sfde) of the form%
\begin{equation}
\left\{
\begin{array}
[c]{ll}%
dX(t) & =b(t,X(t),X_{t},M(t),M_{t},u(t),u_{t})dt+\sigma(t,X(t),X_{t}%
,M(t),M_{t},u(t),u_{t})dB(t)\\
& +%
{\textstyle\int_{\mathbb{R}_{0}}}
\gamma(t,X(t),X_{t},M(t),M_{t},u(t),u_{t},\zeta)\widetilde{N}(dt,d\zeta
);t\in\left[  0,T\right]  ,\\
X(t) & =\xi(t);t\in\left[  -\delta,0\right]  ,\\
u(t) & =u_{0}(t);t\in\left[  -\delta,0\right]  ,
\end{array}
\right.  \label{mfsfde}%
\end{equation}
on a filtered probability space $(%
\Omega
,\mathbb{F},\mathbb{P})$ satisfying the usual conditions, i.e. the filtration
$\mathbb{F}=(\mathcal{F}_{t})_{t\geq0}$ is right-continuous and increasing,
and each $\mathcal{F}_{t}$, $t\geq0$, contains all $\mathbb{P}$-null sets in
$\mathbb{F}$. Here $M(t):=\mathcal{L}(X(t))$ is the law of $X(t)$ at time $t$,
$\delta\geq0$ is a given (constant) memory span and
\begin{equation}
X_{t}:=\{X(t-s)\}_{s\in [0,\delta]}%
\end{equation}
is the path segment of the state process $X(\cdot)$, while
\begin{equation}
M_{t}:=\{M(t+s)\}_{s\in [0,\delta]}%
\end{equation}
is the path segment of the law process $M(\cdot)=\mathcal{L}(X(\cdot))$. The process $u(t)$ is our control process, and $u_t:=\{u(t-s)\}_{s\in[0,\delta]}$ is its memory path segment.
The path processes $X_t, M_t$ and $u_t$ represent the memory terms of the equation \eqref{mfsfde}. The terms
$B(t)$ and $\tilde{N}(dt,d\zeta)$ in the mf-sfde $\left(  \ref{mfsfde}\right)
$ denote a one-dimensional Brownian motion and an independent compensated
Poisson random measure, respectively, such that%
\[%
\begin{array}
[c]{lll}%
\tilde{N}(dt,d\zeta) & = & N(dt,d\zeta)-\nu(d\zeta)dt
\end{array}
\]
where $N(dt,d\zeta)$ is an independent Poisson random measure and $\nu
(d\zeta)$ is the Lévy measure of $N$. For the sake of simplicity, we only
consider the one-dimensional case, i.e. $X(t)\in%
\mathbb{R}
,B(t)\in%
\mathbb{R}
$ and $N(t,\zeta)\in%
\mathbb{R}
,$ for all $t,\zeta$.
\vskip 0.3cm
Let $\mathcal{S}_{\bar{x}}=\mathcal{S}_{\bar{x}}[0,\delta]=\mathbb{R}^{[0,\delta]}$ denote the space
of functions $\bar{x}:[0,\delta]\mapsto
\mathbb{R}$ such that
\[
||\bar{x}||_{\mathcal{S}_{\bar{x}}}^{2}:=
\int_{0}^{\delta}x^{2}(s)ds<\infty.
\]
The spaces $\mathbb{R}^{[0,T]}; T>0$ and $\mathbb{R}^{[0,\infty)}$ are defined similarly.

\begin{definition}{[Segments of elements of $\mathbb{R}^{[\delta,\infty)}$]}
\begin{itemize}
\item 
If $\bar{x} \in \mathbb{R}^{[-\delta,\infty)}$ and $t>0$, we define its \emph{backward/memory path} $\bar{x}_t \in \mathbb{R}^{[0,\delta]}$ by
\begin{equation}
\bar{x}_t (s)=\bar{x}(t-s);\quad s \in [0,\delta].
\end{equation}
\item
If $\bar{x} \in \mathbb{R}^{[-\delta,\infty)}$ and $t>0$, we define its \emph{forward path} $\bar{x}^t \in \mathbb{R}^{[0,\delta]}$ by
\begin{equation}
\bar{x}^t (s)=\bar{x}(t+s);\quad s \in [0,\delta].
\end{equation}
\end{itemize}
\end{definition}

Following Agram and Øksendal \cite{AO1}, we now introduce the following
Hilbert spaces:

\begin{definition}
\-

\begin{itemize}
\item $\mathcal{M}$ is the pre-Hilbert space of random measures $\mu$ on
$\mathbb{R}$ equipped with the norm
\[%
\begin{array}
[c]{lll}%
\left\Vert \mu\right\Vert _{\mathcal{M}}^{2} & := & \mathbb{E[}%
{\textstyle\int_{\mathbb{R}}}
|\hat{\mu}(y)|^{2}e^{-y^{2}}dy]\text{,}%
\end{array}
\]
where $\hat{\mu}$ is the Fourier transform of the measure $\mu$, i.e.%
\[%
\begin{array}
[c]{lll}%
\hat{\mu}(y) & := & {%
{\textstyle\int_{\mathbb{R}}}
}e^{-ixy}d\mu(x);\quad y\in\mathbb{R}.
\end{array}
\]

\item $\mathcal{M}^{\delta}$ is the pre-Hilbert space of all path segments
$\overline{\mu}=\{\mu(s)\}_{s\in\ [0,\delta]}$ of processes $\mu(\cdot)$
with $\mu(s)\in\mathcal{M}$ for each $s\in [0,\delta]$, equipped with
the norm
\begin{equation}
\left\Vert \overline{\mu}\right\Vert ^2_{\mathcal{M}^{\delta}}:={%
{\textstyle\int_{0}^{\delta}}
}\left\Vert \mu(s)\right\Vert ^2_{\mathcal{M}}ds.
\end{equation}

\item $\mathcal{M}_{0}$ and $\mathcal{M}_{0}^{\delta}$ denote the set of
deterministic elements of $\mathcal{M}$ and $\mathcal{M}^{\delta}$,
respectively.\newline For simplicity of notation, in some contexts we regard
$\mathcal{M}$ as a subset of $\mathcal{M}^{\delta}$ and $\mathcal{M}_{0}$ as a
subset of $\mathcal{M}^{\delta}$.
\end{itemize}
\end{definition}

The structure of this space $\mathcal{M}$ equipped with the norm obtained by
the Fourier transform, is an alternative to the Wasserstein metric space
$\mathcal{P}_{2}$ equipped with the Wasserstein distance $W_{2}$. Moreover,
the pre-Hilbert space $\mathcal{M}$ deals with any random measure on $\mathbb{R}$,
however the Wasserstein space $\mathcal{P}_{2}$ deals with Borel probability
measures on $\mathbb{R}$ with finite second moments.

Using a Hilbert space structure for this type of problems has been proposed
by P.L. Lions, to simplify the technicalities of the Wasserstein metric space
where he considers the Hilbert space of square integrable random variables.
Our pre-Hilbert space, however, is new.

In the following, we let $\mathcal{C}:=\mathbb{R}^{[0,\delta]}$ denote the Banach space of all real valued paths 
$\bar{x}:=\{x(s)\}_{s\in\left[  0,\delta \right]  }
$, equipped with the norm%
\begin{equation}
||\bar{x}|| _{\mathcal{C}}:=\mathbb{E}[\underset{s\in\left[
0,\delta \right]  }{\sup}\left\vert x(s)\right\vert ].
\end{equation}
To simplify the writing, we introduce some notations. The
coefficients are assumed to have the form
\[%
\begin{array}
[c]{lll}%
b(t,x,\overline{x},m,\overline{m},u,\overline{u}) & =b(t,x,\overline
{x},m,\overline{m},u,\overline{u},\omega) & :E\rightarrow%
\mathbb{R}
,\\
\sigma(t,x,\overline{x},m,\overline{m},u,\overline{u}) & =\sigma
(t,x,\overline{x},m,\overline{m},u,\overline{u},\omega) & :E\rightarrow%
\mathbb{R}
,\\
\gamma(t,x,\overline{x},m,\overline{m},u,\overline{u},\zeta) & =\gamma
(t,x,\overline{x},m,\overline{m},u,\overline{u},\zeta,\omega) & :E^{\prime
}\rightarrow%
\mathbb{R}
,
\end{array}
\]
where $E:=\left[  0,T\right]  \times%
\mathbb{R}
\times\mathcal{C}\times\mathcal{M}_{0}\times\mathcal{M}_{0}^{\delta}%
\times\mathbb{R}\times\mathcal{C}\times\Omega$ and $E^{\prime}:=\left[
0,T\right]  \times%
\mathbb{R}
\times\mathcal{C}\times\mathcal{M}_{0}\times\mathcal{M}_{0}^{\delta}%
\times\mathbb{R}\times\mathcal{C}\times%
\mathbb{R}
_{0}\times\Omega$ and $\mathbb{R}_{0}=\mathbb{R}\setminus\{0\}$. \newline\newline We
remark that the functionals $b,\sigma$ and $\gamma$ on the mf-sfde depend not just of the solution $X(t)$ and its law $M(t)=\mathcal{L}(X(t))$, but also on the
segment $X_{t}$ and the law of this segment $\mathcal{L}(X_{t})$. This is a
new-type of mean-field stochastic functional differential equations with
memory.\newline\newline 
Let us give some examples: Let $X(t)$ satisfies the
following mean-field delayed sfde%
\begin{equation}
\left\{
\begin{array}
[c]{ll}%
dX(t) & =b(t,\mathbf{X}(t),\mathbb{E}[\mathbf{X}(t)],u(t))dt+\sigma
(t,\mathbf{X}(t),\mathbb{E}[\mathbf{X}(t)],u(t))dB(t)\\
& +%
{\textstyle\int_{\mathbb{R}_{0}}}
\gamma(t,\mathbf{X}(t),\mathbb{E}[\mathbf{X}(t)],u(t),\zeta)\tilde
{N}(dt,d\zeta);t\in\left[  0,T\right]  ,\\
X(t) & =\xi(t);t\in\left[  -\delta,0\right]  ,
\end{array}
\right.  \label{mfd}%
\end{equation}
where we denote by the bold $\mathbf{X}(t)=%
{\textstyle\int_{0}^{\delta}}
X(t-s)\mu(ds)$ for some bounded Borel-measure $\mu$. As noted in Agram and
Røse \cite{AR} and Banos et al \cite{G et al}, we have the following:

\begin{itemize}
\item If this measure $\mu$ is a Dirac-measure concentrated at $0$ i.e.
$\mathbf{X}(t)=X(t)$ then equation $\left(  \ref{mfd}\right)  $ is a
\emph{classical mean-field stochastic differential equation}, we refer for
example to Anderson and Djehiche in \cite{AD} and Hu el al in \cite{HOS} for
stochastic control of such a systems.

\item It could also be the Dirac measure concentrated at $\delta$
then\ $\mathbf{X}(t)=X(t-\delta)$ and in that case the state equation is
called a \emph{mean-field sde with discrete delay}, see for instance Meng and
Shen \cite{MS} and for delayed systems without a mean-field term, we refer to
Chen and Wu \cite{CW}, Dahl et al \cite{DMOR} and Øksendal et al \cite{osz}.

\item If we choose now $\mu(ds)=g(s)ds$ for any function $g\in L^{1}(\left[
0,\delta\right]  )$ thus $\mathbf{X}(t)=%
{\textstyle\int_{0}^{\delta}}
g(s)X(t-s)ds$ and the state is a \emph{mean-field distributed delay}.\newline
\end{itemize}

It is worth mentioning the papers by Lions \cite{Lions}, Cardaliaguet
\cite{Cardaliaguet}, Carmona and Delarue \cite{carmona1}, \cite{carmona},
Buckdahn et al \cite{BLPR} and Agram \cite{A} for more details about systems
driven by mean-field equations and stochastic control problems for such a
system. These papers, however, use the Wasserstein metric space of probability
measures and not our Hilbert space of measures.\newline\newline

The paper is organized as follows: In section 2, we give some mathematical
background and define some concepts and spaces which will be used in the
paper. In section 3, we prove existence and uniqueness of memory McKean-Vlasov
equations. Section 4 contains the main results of this paper, including a
sufficient and a necessary maximum principle for the optimal control of
stochastic memory mean-field equations. In section 5, we illustrate our
results by solving a mean-variance and a linear-quadratic problems of a memory processes.

\section{Generalities}

In this section, we recall some concepts which will be used on the sequel.

\begin{description}
\item[a)] We first discuss the differentiability of functions defined on a
Banach space.
\end{description}

Let $\mathcal{X},\mathcal{Y}$ be two Banach spaces with norms $\Vert\cdot
\Vert_{\mathcal{X}},\Vert\cdot\Vert_{\mathcal{Y}}$, respectively, and let
$F:\mathcal{X}\rightarrow\mathcal{Y}$.

\begin{itemize}
\item We say that $F$ has a directional derivative (or Gâteaux derivative) at
$v\in\mathcal{X}$ in the direction $w\in\mathcal{X}$ if
\[
D_{w}F(v):=\lim_{\varepsilon\rightarrow0}\frac{1}{\varepsilon}(F(v+\varepsilon
w)-F(v))
\]
exists.

\item We say that $F$ is Fréchet differentiable at $v\in\mathcal{X}$ if there
exists a continuous linear map $A:\mathcal{X}\rightarrow\mathcal{Y}$ such
that
\[
\lim_{\substack{h\rightarrow0\\h\in\mathcal{X}}}\frac{1}{\Vert h\Vert
_{\mathcal{X}}}\Vert F(v+h)-F(v)-A(h)\Vert_{\mathcal{Y}}=0,
\]
where $A(h)=\langle A,h\rangle$ is the action of the linear operator $A$ on
$h$. In this case we call $A$ the \textit{gradient} (or Fréchet derivative) of
$F$ at $v$ and we write
\[
A=\nabla_{v}F.
\]

\item If $F$ is Fréchet differentiable at $v$ with Fréchet derivative
$\nabla_{v}F$, then $F$ has a directional derivative in all directions
$w\in\mathcal{X}$ and
\[
D_{w}F(v)=\nabla_{v}F(w)=\langle\nabla_{v}F,w\rangle.
\]
In particular, note that if $F$ is a linear operator, then $\nabla_{v}F=F$ for
all $v$.\newline
\end{itemize}

\begin{description}
\item[b)] Throughout this work, we will use the following spaces:
\end{description}

\begin{itemize}
\item $\mathcal{S}^{2}$ is the set of ${\mathbb{R}}$-valued $\mathbb{F}%
$-adapted càdlàg processes $(X(t))_{t\in\lbrack-\delta,T]}$ such that
\[
{\Vert X\Vert}_{\mathcal{S}^{2}}^{2}:={\mathbb{E}}[\sup_{t\in\lbrack
-\delta,T]}|X(t)|^{2}]~<~\infty\;,
\]
(alternatively $(X(t))_{t\in\lbrack0,T+\delta]}$ with
\[
{\Vert X\Vert}_{\mathcal{S}^{2}}^{2}={\mathbb{E}}[\sup_{t\in\lbrack
0,T+\delta]}|X(t)|^{2}]~<~\infty\;,
\]
depending on the context.)

\item $\mathbb{L}^{2}$ is the set of ${\mathbb{R}}$-valued $\mathbb{F}%
$-adapted processes $(Q(t))_{t\in\lbrack0,T]}$ such that
\[
\Vert Q\Vert_{\mathbb{L}^{2}}^{2}:={\mathbb{E}}[%
{\textstyle\int_{0}^{T}}
|Q(t)|^{2}dt]<~\infty\;.
\]

\item $\mathcal{U}^{ad}$ is a set of all stochastic processes $u$ required to
have values in a convex subset $\mathcal{U}$ of $%
\mathbb{R}
$ and adapted to a given subfiltration $\mathbb{G}=\{\mathcal{G}_{t}%
\}_{t\geq0},$ where $\mathcal{G}_{t}\subseteq\mathcal{F}_{t}$ for all $t\geq
0$. We call $\mathcal{U}^{ad}$ the set of admissible control processes
$u(\cdot)$.

\item $L^{2}(\mathcal{F}_{t})$ is the set of ${\mathbb{R}}$-valued square
integrable $\mathcal{F}_{t}$-measurable random variables.

\item $\mathbb{L}_{\nu}^{2}$ is the set of ${\mathbb{R}}$-valued $\mathbb{F}%
$-adapted processes $Z:\mathbb{R}_{0}\rightarrow%
\mathbb{R}
$ such that
\[
||Z||_{\mathbb{L}_{\nu}^{2}}^{2}:={\mathbb{E}}[%
{\textstyle\int_{\mathbb{R}_{0}}}
|Z(t,\zeta)|^{2}\nu(d\zeta)dt]~<~\infty\;.
\]

\item $\mathcal{R}$ is the set of measurable functions $r:\mathbb{R}%
_{0}\rightarrow\mathbb{R}.$

\item $C_{a}([0,T],\mathcal{M}_{0})$ denotes the set of absolutely continuous
functions $m:[0,T]\rightarrow\mathcal{M}_{0}.$

\end{itemize}

\section{Solvability of memory mean-field sfde}

For a given constant $\delta>0$, we consider a memory mean-field stochastic
functional differential equations (mf-sfde) of the following form:
\begin{equation}
\left\{
\begin{array}
[c]{ll}%
dX(t) & =b(t,X(t),X_{t},M(t),M_{t})dt+\sigma(t,X(t),X_{t},M(t),M_{t})dB(t)\\
& +%
{\textstyle\int_{\mathbb{R}_{0}}}
\gamma(t,X(t),X_{t},M(t),M_{t},\zeta)\widetilde{N}(dt,d\zeta);t\in\left[
0,T\right]  ,\\
X(t) & =\xi(t);t\in\left[  -\delta,0\right]  .
\end{array}
\right.  \label{sfde}%
\end{equation}
Here $E:=\left[  0,T\right]  \times%
\mathbb{R}
\times\mathcal{C}\times\mathcal{M}_{0}\times\mathcal{M}_{0}^{\delta}%
\times\Omega$, $E^{\prime}:=\left[  0,T\right]  \times%
\mathbb{R}
\times\mathcal{C}\times\mathcal{M}_{0}\times\mathcal{M}_{0}^{\delta}\times%
\mathbb{R}
_{0}\times\Omega$ and the coefficients
\[%
\begin{array}
[c]{lll}%
b(t,x,\overline{x},m,\overline{m}) & =b(t,x,\overline{x},m,\overline{m}%
,\omega) & :E\rightarrow%
\mathbb{R}
,\\
\sigma(t,x,\overline{x},m,\overline{m}) & =\sigma(t,x,\overline{x}%
,m,\overline{m},\omega) & :E\rightarrow%
\mathbb{R}
,\\
\gamma(t,x,\overline{x},m,\overline{m},\zeta) & =\gamma(t,x,\overline
{x},m,\overline{m},\zeta,\omega) & :E^{\prime}\rightarrow%
\mathbb{R}
,
\end{array}
\]
are supposed to be $\mathcal{F}_{t}$-measurable and the initial value function
$\xi$ is $\mathcal{F}_{0}$-measurable.\newline

For more information about stochastic functional differential equations, we
refer to the seminal work of S.E.A. Mohammed \cite{M} and a recent paper by
Banos et al \cite{G et al}.

In order to prove an existence and uniqueness result for the mf-sfde $\left(
\ref{sfde}\right)  $, we first need the following lemma:

\begin{lemma}
\-

\begin{description}
\item[(i)] Let $X^{(1)}$ and $X^{(2)}$ be two random variables in
$L^{2}(\mathbb{P})$. Then
\[%
\begin{array}
[c]{lll}%
\left\Vert \mathcal{L}(X^{(1)})-\mathcal{L}(X^{(2)})\right\Vert _{\mathcal{M}%
_{0}}^{2} & \leq & \sqrt{\pi}\mathbb{E}[(X^{(1)}-X^{(2)})^{2}]\text{.}%
\end{array}
\]

\item[(ii)] Let $\{X^{(1)}(t)\}_{t\geq0},$ $\{X^{(2)}(t)\}_{t\geq0}$ be two
processes such that
\[
\mathbb{E}[%
{\textstyle\int_{0}^{T}}
X^{(i)2}(s)ds]<\infty\text{ for all }T\text{ with }i=1,2\text{.}%
\]
Then, for all $t$,
\[%
\begin{array}
[c]{lll}%
||\mathcal{L}(X_{t}^{(1)})-\mathcal{L}(X_{t}^{(2)})||_{\mathcal{M}_{0}%
^{\delta}}^{2} & \leq & \sqrt{\pi}\mathbb{E}[%
{\textstyle\int_{0}^{\delta}}
(X^{(1)}(t-s)-X^{(2)}(t-s))^{2}ds]\text{.}%
\end{array}
\]

\end{description}
\end{lemma}

\noindent{Proof.} \quad By definition of the norms and standard properties of
the complex exponential function, we have%
\[%
\begin{array}
[c]{l}%
\left\Vert \mathcal{L}(X^{(1)})-\mathcal{L}(X^{(2)})\right\Vert _{\mathcal{M}%
_{0}}^{2}\\
=%
{\textstyle\int_{\mathbb{R}}}
|\widehat{\mathcal{L}}(X^{(1)})(y)-\widehat{\mathcal{L}}(X^{(2)}%
)(y)|^{2}e^{-y^{2}}dy\\
=%
{\textstyle\int_{\mathbb{R}}}
\left\vert
{\textstyle\int_{\mathbb{R}}}
e^{-ixy}d\mathcal{L}(X^{(1)})(x)-%
{\textstyle\int_{\mathbb{R}}}
e^{-ixy}d\mathcal{L}(X^{(2)})(x)\right\vert ^{2}e^{-y^{2}}dy\\
=%
{\textstyle\int_{\mathbb{R}}}
|\mathbb{E}[e^{-iyX^{(1)}}-e^{-iyX^{(2)}}]|^{2}e^{-y^{2}}dy\\
=%
{\textstyle\int_{\mathbb{R}}}
|\mathbb{E}[\cos(yX^{(1)})-\cos(yX^{(2)})]-i\mathbb{E}[\sin(yX^{(1)}%
)-\sin(yX^{(2)})]|^{2}e^{-y^{2}}dy\\
=%
{\textstyle\int_{\mathbb{R}}}
(\mathbb{E}[\cos(yX^{(1)})-\cos(yX^{(2)})]^{2}+\mathbb{E}[\sin(yX^{(1)}%
)-\sin(yX^{(2)})]^{2})e^{-y^{2}}dy\\
\leq%
{\textstyle\int_{\mathbb{R}}}
(\mathbb{E}[|\cos(yX^{(1)})-\cos(yX^{(2)})|]^{2}+\mathbb{E}[|\sin
(yX^{(1)})-\sin(yX^{(2)})|]^{2})e^{-y^{2}}dy\\
\leq%
{\textstyle\int_{\mathbb{R}}}
(\mathbb{E}[|y(X^{(1)}-X^{(2)})|]^{2}+\mathbb{E}[|y(X^{(1)}-X^{(2)}%
)|]^{2})e^{-y^{2}}dy\\
\leq2%
{\textstyle\int_{\mathbb{R}}}
y^{2}e^{-y^{2}}dy\mathbb{E}[|X^{(1)}-X^{(2)}|]^{2}\\
\leq\sqrt{\pi}\mathbb{E}[(X^{(1)}-X^{(2)})^{2}].
\end{array}
\]

Similarly we get that
\[%
\begin{array}
[c]{lll}%
||\mathcal{L}(X_{t}^{(1)})-\mathcal{L}(X_{t}^{(2)})||_{\mathcal{M}_{0}%
^{\delta}}^{2} & \leq &
{\textstyle\int}_0^{\delta}
\left\Vert \mathcal{L}(X^{(1)}(t-s))-\mathcal{L}(X^{(2)}(t-s))\right\Vert _{\mathcal{M}_{0}%
}^{2}ds\\
& \leq & \sqrt{\pi}\mathbb{E}[%
{\textstyle\int}_0^{\delta}
(X^{(1)}(t-s)-X^{(2)}(t-s))^{2}ds].
\end{array}
\]
$\square$\\

We also need the following result, which is Lemma 2.3 in \cite{AO1}:
\begin{lemma}
Suppose that $X(t)$ is an Itô-Lévy process of the form
\begin{equation}%
\begin{cases}
dX(t)=\alpha(t)dt+\beta(t)dB(t)+{\textstyle\int_{\mathbb{R}_{0}}}
\gamma(t,\zeta)\tilde{N}(dt,d\zeta);\quad t\in\lbrack0,T],\\
X(0)=x\in\mathbb{R},
\end{cases}
\label{eq2.1}%
\end{equation}
where $\alpha,\beta$ and $\gamma$ are predictable processes.\\
Then the map $t\mapsto M(t):[0,T]\rightarrow\mathcal{M}_{0}$ is absolutely continuous.
\end{lemma}
\vskip 0.3cm
It follows that $t \mapsto M(t)$ is differentiable for a.a.$t$. We will in the following use the notation
\begin{equation}
M'(t)=\frac{dM(t)}{dt}.
\end{equation}
We are now able to state the theorem of existence and uniqueness of a solution
of equation $\left(  \ref{sfde}\right)  $. As before we put $E:=\left[
0,T\right]  \times%
\mathbb{R}
\times\mathcal{C}\times\mathcal{M}_{0}\times\mathcal{M}_{0}^{\delta}%
\times\Omega$ and $E^{\prime}:=\left[  0,T\right]  \times%
\mathbb{R}
\times\mathcal{C}\times\mathcal{M}_{0}\times\mathcal{M}_{0}^{\delta}\times%
\mathbb{R}
_{0}\times\Omega$. Then we have

\begin{theorem}
\label{existence}Assume that $\xi(t)\in\mathcal{C},b,\sigma:E\rightarrow%
\mathbb{R}
$ and $\gamma:E^{\prime}\rightarrow%
\mathbb{R}
$ are progressively measurable and satisfy the following uniform Lipschitz
condition $dtP(d\omega)$-a.e.:\newline There is some constant $L\in%
\mathbb{R}
$ such that%
\begin{equation}%
\begin{array}
[c]{l}%
|b(t,x,\overline{x},m,\overline{m},\omega)-b(t,x^{\prime},\overline{x}%
^{\prime},m^{\prime},\overline{m}^{\prime},\omega)|^{2}+|\sigma(t,x,\overline
{x},m,\overline{m},\omega)-\sigma(t,x^{\prime},\overline{x}^{\prime}%
,m^{\prime},\overline{m}^{\prime},\omega)|^{2}\\
\text{ \ \ \ \ \ \ \ }+\int_{%
\mathbb{R}
_{0}}|\gamma(t,x,\overline{x},m,\overline{m},\zeta,\omega)-\gamma(t,x^{\prime
},\overline{x}^{\prime},m^{\prime},\overline{m}^{\prime},\zeta,\omega)|^{2}%
\nu(d\zeta)\\
\leq L(|x-x^{\prime}|^{2}+||\overline{x}-\overline{x}^{\prime}||_{\mathcal{C}%
}^{2}+||m-m^{\prime}||_{\mathcal{M}_{0}}^{2}+||\overline{m}-\overline
{m}^{\prime}||_{\mathcal{M}_{0}^{\delta}}^{2}),\text{ for a.a. }t,\omega,
\end{array}
\label{Lip}%
\end{equation}
and%
\begin{equation}%
\begin{array}
[c]{l}%
|b(t,0,0,\mu_{0},\mu_{0},\omega)|^{2}+|\sigma(t,0,0,\mu_{0},\mu_{0}%
,\omega)|^{2}\\
\text{ \ }+\int_{%
\mathbb{R}
_{0}}|\gamma(t,0,0,\mu_{0},\mu_{0},\zeta,\omega)|^{2}\nu(d\zeta)\leq L\text{
for a.a. }t,\omega,
\end{array}
\label{Bou}%
\end{equation}
where $\mu_{0}$ is the Dirac measure with mass at zero. Then there is a unique
solution $X\in\mathcal{S}^{2}$ of the mf-sfde $\left(  \ref{sfde}\right)  $.
\end{theorem}

\noindent{Proof.} \quad For $X\in\mathcal{S}^{2}[-\delta,T]$ and for $t_{0}%
\in(0,T]$, we introduce the norm%
\[
||X||_{t_{0}}^{2}:=\mathbb{E[}\underset{t\in\lbrack-\delta,t_{0}]}{\sup
}|X(t)|^{2}].
\]

The space $\mathbb{H}_{t_{0}}$ equipped with this norm is a Banach space.
Define the mapping $\Phi:\mathbb{H}_{t_{0}}\rightarrow\mathbb{H}_{t_{0}}$ by
$\Phi(x)=X$ where $X\in\mathcal{S}^{2}$ is defined by%
\[
\left\{
\begin{array}
[c]{ll}%
dX(t) & =b(t,x(t),x_{t},m(t),m_{t})dt+\sigma(t,x(t),x_{t},m(t),m_{t})dB(t)\\
& +%
{\textstyle\int_{\mathbb{R}_{0}}}
\gamma(t,x(t),x_{t},m(t),m_{t},\zeta)\widetilde{N}(dt,d\zeta);t\in\left[
0,T\right]  ,\\
X(t) & =\xi(t);t\in\left[  -\delta,0\right]  .
\end{array}
\right.
\]
We want prove that $\Phi$ is contracting in $\mathbb{H}_{t_{0}}$ under the
norm $||\cdot||_{t_{0}}$ for small enough $t_{0}$. For two arbitrary elements
$(x^{1},x^{2})$ and $(X^{1},X^{2})$, we denote their difference by
$\widetilde{x}=x^{1}-x^{2}$ and $\widetilde{X}=X^{1}-X^{2}$ respectively. In
the following $C<\infty$ will denote a constant which is big enough for all
the inequalities to hold.

Applying the Itô formula to $\widetilde{X}^{2}(t)$, we get%

\begin{align*}
\widetilde{X}^{2}(t) &  =2%
{\textstyle\int_{0}^{t}}
\widetilde{X}(s)(b(s,x^{1}(s),x_{s}^{1},m^{1}(s),m_{s}^{1})-b(s,x^{2}%
(s),x_{s}^{2},m^{2}(s),m_{s}^{2}))ds\\
&  +2%
{\textstyle\int_{0}^{t}}
\widetilde{X}(s)(\sigma(s,x^{1}(s),x_{s}^{1},m^{1}(s),m_{s}^{1})-\sigma
(s,x^{2}(s),x_{s}^{2},m^{2}(s),m_{s}^{2}))dB(s)\\
&  +2%
{\textstyle\int_{0}^{t}}
\widetilde{X}(s)%
{\textstyle\int_{\mathbb{R}_{0}}}
(\gamma(s,x^{1}(s),x_{s}^{1},m^{1}(s),m_{s}^{1},\zeta)-\gamma(s,x^{2}%
(s),x_{s}^{2},m^{2}(s),m_{s}^{2},\zeta))\widetilde{N}(ds,d\zeta)\\
&  +%
{\textstyle\int_{0}^{t}}
(\sigma(s,x^{1}(s),x_{s}^{1},m^{1}(s),m_{s}^{1})-\sigma(s,x^{2}(s),x_{s}%
^{2},m^{2}(s),m_{s}^{2}))^{2}ds\\
&  +%
{\textstyle\int_{0}^{t}}
{\textstyle\int_{\mathbb{R}_{0}}}
(\gamma(s,x^{1}(s),x_{s}^{1},m^{1}(s),m_{s}^{1},\zeta)-\gamma(s,x^{2}%
(s),x_{s}^{2},m^{2}(s),m_{s}^{2},\zeta))^{2}\nu(d\zeta)ds.
\end{align*}
By the Lipschitz assumption $\left(  \ref{Lip}\right)  $ combined with
standard majorization of the square of a sum (resp. integral) via the sum
(resp. integral) of the square (up to a constant), we get%

\begin{align*}
\widetilde{X}^{2}(t)  &  \leq C%
{\textstyle\int_{0}^{t}}
|\widetilde{X}(s)|\Delta_{t_{0}}ds\\
&  +|%
{\textstyle\int_{0}^{t}}
\widetilde{X}(s)\widetilde{\sigma}(s)dB(s)|+|%
{\textstyle\int_{0}^{t}}
{\textstyle\int_{\mathbb{R}_{0}}}
\widetilde{X}(s)\widetilde{\gamma}(s,\zeta)\widetilde{N}(ds,d\zeta
)|+tC\Delta_{t_{0}}^{(2)},
\end{align*}
where
\[%
\begin{array}
[c]{ll}%
\Delta_{t_{0}} & :=||\widetilde{x}||_{\mathcal{S}^{2}}+||\widetilde
{\overline{x}}||_{\mathcal{C}}+||\widetilde{m}||_{\mathcal{M}_{0}%
}+||\widetilde{\overline{m}}||_{\mathcal{M}_{0}^{\delta}}\\
\Delta_{t_{0}}^{(2)} & :=||\widetilde{x}||_{\mathcal{S}^{2}}^{2}%
+||\widetilde{\overline{x}}||_{\mathcal{C}}^{2}+||\widetilde{m}||_{\mathcal{M}%
_{0}}^{2}+||\widetilde{\overline{m}}||_{\mathcal{M}_{0}^{\delta}}^{2}.
\end{array}
\]
By the Burkholder-Davis-Gundy inequalities,
\begin{equation}
\mathbb{E}[\sup_{t\leq t_{0}}|%
{\textstyle\int_{0}^{t}}
\widetilde{X}(s)\widetilde{\sigma}(s)dB(s)|]\leq C\mathbb{E}[(%
{\textstyle\int_{0}^{t_{0}}}
\widetilde{X}^{2}(s)\widetilde{\sigma}^{2}(s)ds)^{\frac{1}{2}}]\leq
Ct_{0}||\widetilde{X}||_{t_{0}}\Delta_{t_{0}},
\end{equation}
and
\begin{equation}
\mathbb{E}[\sup_{t\leq t_{0}}|%
{\textstyle\int_{0}^{t}}
\widetilde{X}(s)\widetilde{\gamma}(s)\widetilde{N}(ds,d\zeta)|]\leq
C\mathbb{E}[(%
{\textstyle\int_{0}^{t_{0}}}
\widetilde{X}^{2}(s)\widetilde{\gamma}^{2}(s)\nu(d\zeta)ds)^{\frac{1}{2}}]\leq
Ct_{0}||\widetilde{X}||_{t_{0}}\Delta_{t_{0}}.
\end{equation}
Combining the above and using that
\[
||\widetilde{X}||_{t_{0}}\Delta_{t_{0}}\leq C(||X||_{t_{0}}^{2}+\Delta_{t_{0}%
}^{(2)}),
\]
we obtain%

\[
||\widetilde{X}||_{t_{0}}^{2}:=\mathbb{E}[\sup_{t\leq t_{0}}\widetilde{X}%
^{2}(t)]\leq Ct_{0}(||\widetilde{X}||_{t_{0}}^{2}+\Delta_{t_{0}}^{(2)}).
\]
By definition of the norms, we have
\begin{equation}
\Delta_{t_{0}}^{(2)}\leq C||\widetilde{x}||_{t_{0}}^{2}.
\end{equation}
Thus we see that if $t_{0}>0$ is small enough we obtain
\begin{equation}
||\widetilde{X}(t)||_{t_{0}}^{2}\leq\frac{1}{2}||\widetilde{x}(s)||_{t_{0}%
}^{2},
\end{equation}
and hence $\Phi$ is a contraction on $\mathbb{H}_{t_{0}}$. Therefore the
equation has a solution up to $t_{0}$. By the same argument we see that the
solution is unique. Now we repeat the argument above, but starting at $t_{0}$
instead of starting at $0$. Then we get a unique solution up to $2t_{0}$.
Iterating this, we obtain a unique solution up to $T$ for any $T<\infty.$
\hfill$\square$

\section{Optimal control of memory mf-sfde}

Consider again the controlled memory mf-sfde $\left(  \ref{mfsfde}\right)  $
\begin{equation}
\left\{
\begin{array}
[c]{ll}%
dX(t) & =b(t,X(t),X_{t},M(t),M_{t},u(t),u_t)dt+\sigma(t,X(t),X_{t},M(t),M_{t},u(t),u_t)dB(t)\\
& +%
{\textstyle\int_{\mathbb{R}_{0}}}
\gamma(t,X(t),X_{t},M(t),M_{t},u(t),u_t,\zeta)\widetilde{N}(dt,d\zeta);t\in\left[
0,T\right]  ,\\
X(t) & =\xi(t);t\in\left[  -\delta,0\right]  .
\end{array}
\right.  \label{exmfsfde}%
\end{equation}

The coefficients $b,\sigma$ and $\gamma$ are supposed to satisfy the
assumptions of Theorem \ref{existence}, uniformly w.r.t. $u\in\mathcal{U}^{ad}$,
then we have the existence and the uniqueness of the solution $X(t)\in
\mathcal{S}^{2}$ of the controlled mf-sfde $\left(  \ref{exmfsfde}\right)  $.

Moreover, $b,\sigma$ and $\gamma$ have Fréchet derivatives w.r.t.
$\overline{x}$, $m$, $\overline{m}$ and are continuously differentiable in the
variables $x$ and $u$.\newline

The performance functional is assumed to be of the form%
\begin{equation}%
\begin{array}
[c]{lll}%
J(u) & = & \mathbb{E[}%
{\textstyle\int_{0}^{T}}
\ell(t,X(t),X_{t},M(t),M_{t},u(t),u_{t})dt+h(X(T),M(T)]\text{; }%
u\in\mathcal{U}.
\end{array}
\label{perf}%
\end{equation}
With $E:=\left[  0,T\right]  \times%
\mathbb{R}
\times\mathcal{C}\times\mathcal{M}_{0}\times\mathcal{M}_{0}^{\delta}%
\times\mathcal{U}^{ad}\times\mathcal{C}\times\Omega,$ $E^{\prime}:=%
\mathbb{R}
\times\mathcal{M}_{0}\times\Omega$ we assume that the functions%
\[%
\begin{array}
[c]{lll}%
\ell(t,x,\bar{x},m,\bar{m},u,\bar{u}) & =\ell(t,x,\bar{x},m,\bar{m},u,\bar
{u},\omega) & :E\rightarrow%
\mathbb{R}
,\\
h(x,m) & =h(x,m,\omega) & :E^{\prime}\rightarrow%
\mathbb{R}
,
\end{array}
\]
admit Fréchet derivatives w.r.t. $\overline{x}$, $m$, $\overline{m}$ and are
continuously differentiable w.r.t. $x$ and $u$. We allow the integrand in the
performance functional $\left(  \ref{perf}\right)  $ to depend on the path
process $X_{t}$ and also its law process $\mathcal{L}(X_{t})=:M_{t},$ and we
allow the terminal value to depend on the state $X(T)$ and its law $M(T)$.
\newline\newline Consider the following optimal control problem. It may
regarded as a partial information control problem (since $u$ is required to be
$\mathbb{G}$-adapted) but only in the limited sense, since $\mathbb{G}$ does
not depend on the observation.

\begin{problem}
Find $u^{*} \in\mathcal{U}^{ad}$ such that
\begin{equation}
\label{eq4.3}J(u^{*})= \sup_{u \in\mathcal{U}^{ad}} J(u).
\end{equation}

\end{problem}

To study this problem we first introduce its associated Hamiltonian, as follows:

\begin{definition}
The Hamiltonian
\[
H:[0,T+\delta]\times\mathbb{R}\times\mathcal{C}\times\mathcal{M}_{0}%
\times\mathcal{M}_{0}^{\delta}\times\mathcal{U}^{ad}\times\mathcal{C}%
\times\mathbb{R}\times\mathbb{R}\times\mathcal{R}\times\mathbb{K}\times
\Omega\rightarrow\mathbb{R}%
\]
associated to this memory mean-field stochastic control problem \eqref{eq4.3}
is defined by%
\begin{equation}%
\begin{array}
[c]{ll}%
H(t,x,\overline{x},m,\overline{m},u,\overline{u},p^{0},q^{0},r^{0}%
(\cdot),p^{1}) & =H(t,x,\overline{x},m,\overline{m},u,\overline{u},p^{0}%
,q^{0},r^{0}(\cdot),p^{1},\omega)\\
& =\ell(t,x,\overline{x},m,\overline{m},u,\overline{u})+p^{0}b(t,x,\overline
{x},m,\overline{m},u,\overline{u})\\
& \text{ \ \ \ }+q^{0}\sigma(t,x,\overline{x},m,\overline{m},u,\overline{u})\\
& \text{ \ \ \ }+%
{\textstyle\int_{\mathbb{R}_{0}}}
r^{0}(t,\zeta)\gamma\left(  t,\zeta\right)  \nu(d\zeta)+\left\langle
p^{1},m^{\prime}\right\rangle \text{; }t\in\lbrack0,T],
\end{array}
\label{haml}%
\end{equation}
and $H(t,x,\overline{x},m,\overline{m},u,\overline{u},p^{0},q^{0},r^{0}%
(\cdot),p^{1})=0; \quad t>T$.
\end{definition}

The Hamiltonian $H$ is assumed to be continuously differentiable w.r.t. $x,u$
and to admit Fréchet derivatives w.r.t. $\overline{x},m$ and $\overline{m}$.

In the following we let $L_{0}^{2}$ denote the set of measurable stochastic
processes $Y(t)$ on $\mathbb{R}$ such that $Y(t)=0$ for $t<0$ and for $t>T$
and
\begin{equation}%
{\textstyle\int_{0}^{T}}
Y^{2}(t)dt<\infty\quad a.s.
\end{equation}
\newline The map
\[
Y\mapsto%
{\textstyle\int_{0}^{T}}
<\nabla_{\overline{x}}H(t),Y_{t}>dt;\quad Y\in L_{0}^{2}%
\]
is a bounded linear functional on $L_{0}^{2}$. Therefore, by the Riesz
representation theorem there exists a unique process $\Gamma_{\bar{x}}(t)\in L_{0}^{2}$
such that
\begin{equation}%
{\textstyle\int_{0}^{T}}
\Gamma_{\bar{x}}(t)Y(t)dt=%
{\textstyle\int_{0}^{T}}
<\nabla_{\overline{x}}H(t),Y_{t}>dt, \label{eq4.6}%
\end{equation}
for all $Y\in L_{0}^{2}$. Here $<\nabla_{\overline{x}}H(t),Y_{t}>$ denotes the
action of the operator $\nabla_{\overline{x}}H(t)$ to the segment
$Y_{t}=\{Y(t-s)\}_{s\in [0,\delta]}$, where $H(t)$ is a shorthand
notation for
\[
H(t,X(t),X_{t},M(t),M_{t},u(t),u_{t},p^{0}(t),q^{0}(t),r^{0}(t,\cdot
),p^{1}(t),\omega).
\]
As a suggestive notation (see below) for $\Gamma_{\bar{x}}$ we will in the following
write
\begin{equation}
\nabla_{\overline{x}}H^{t}:=\Gamma_{\bar{x}}(t).
\end{equation}

\begin{lemma}
Consider the case when
\[
H(t,x,\overline{x},p,q)= f(t,x) + F(\overline{x})p +\sigma q,
\]
Then
\begin{equation}\label{eq4.8}
\Gamma_{\bar{x}}(t):= <\nabla_{\overline{x}}F,p^{t}>
\end{equation}
satisfies \eqref{eq4.6}, where $p^{t}:=\{p(t+r)\}_{r \in[0,\delta
]}$.
\end{lemma}

\noindent{Proof.} \quad \quad We must verify that if we define $\Gamma_{\bar{x}}(t)$ by \eqref{eq4.8}, then \eqref{eq4.6} holds. To this end, choose $Y\in L_{0}^{2}$ and consider%
\[%
\begin{array}
[c]{ll}%
\int_0^T \Gamma_{\bar{x}}(t)Y(t)dt=
{\textstyle\int_{0}^{T}}
<\nabla_{\bar{x}}F,p^{t}>Y(t)dt & =%
{\textstyle\int_{0}^{T}}
<\nabla_{\bar{x}}F,\{p(t+r)\}_{r\in\lbrack0,\delta]}>Y(t)dt\\
& =%
{\textstyle\int_{0}^{T}}
<\nabla_{\bar{x}}F,\{Y(t)p(t+r)\}_{r\in\lbrack0,\delta]}>dt\\
& =<\nabla_{\bar{x}}F,\{%
{\textstyle\int_{r}^{T+r}}
Y(u-r)p(u)du\}_{r\in\lbrack0,\delta]}>\\
& =<\nabla_{\bar{x}}F,\{%
{\textstyle\int_{0}^{T}}
Y(u-r)p(u)du\}_{r\in\lbrack0,\delta]}\\
& =%
{\textstyle\int_{0}^{T}}
<\nabla_{\bar{x}}F,Y_{u}>p(u)du\\
& =%
{\textstyle\int_{0}^{T}}
<\nabla_{\bar{x}}H(u),Y_{u}>du.
\end{array}
\]
\hfill$\square$ \bigskip

\begin{example}
(i) For example, if $a \in \mathbb{R}^{[0,\delta]}$ is a bounded function and $F(\bar{x})$ is the averaging operator defined by
\begin{equation}
F(\bar{x})=%
{\textstyle\int_{0}^{\delta}}
a(s)x(s)ds
\end{equation}
when $\bar{x}=\{x(s)\}_{s\in [0,\delta]}$, then
\begin{equation}
<\nabla_{\bar{x}}F,p^{t}>=<F,p^{t}>=
{\textstyle\int_{0}^{\delta}}
a(r)p(t+r)dr.
\end{equation}

(ii) Similarly, if $t_0 \in [0,\delta]$ and $G$ is evaluation at $t_0$, i.e.
\begin{equation}
G(\bar{x})=x(t_0)\text{ when }\bar{x}=\{x(s)\}_{s\in [0,\delta]},
\end{equation}
then
\begin{equation}
<\nabla_{\bar{x}}G,p^{t}>=p(t+t_0).
\end{equation}

\end{example}

For $u\in\mathcal{U}^{ad}$ with corresponding solution $X=X^{u}$, define
$p=(p^{0},p^{1}),q=(q^{0},q^{1})$ and $r=(r^{0},r^{1})$ by the following two
adjoint equations:\newline

\begin{itemize}
\item The advanced backward stochastic functional differential equation
(absfde) in the unknown $(p^{0},q^{0},r^{0})\in\mathcal{S}^{2}\times
\mathbb{L}^{2}\times\mathbb{L}_{\nu}^{2}$ is given by
\begin{equation}
\left\{
\begin{array}
[c]{ll}%
dp^{0}(t) & =-[\tfrac{\partial H}{\partial x}(t)+\mathbb{E}(\nabla
_{\overline{x}}H^{t}|\mathcal{F}_{t})]dt+q^{0}(t)dB(t)+\int_{%
\mathbb{R}
_{0}}r^{0}(t,\zeta)\widetilde{N}(dt,d\zeta);t\in\lbrack0,T],\\
p^{0}(t) & =\tfrac{\partial h}{\partial x}(X(T),M(T));t\geq T,\\
q^{0}(t) & =0;t>T,\\
r^{0}(t,\cdot) & =0;t>T.
\end{array}
\right.  \label{p0}
\end{equation}

\item The operator-valued mean-field advanced backward stochastic functional
differential equation (ov-mf-absfde) in the unknown $(p^{1},q^{1},r^{1}%
)\in\mathcal{S}_{\mathbb{K}}^{2}\times\mathbb{L}^{2}_{\mathbb{K}}%
\times\mathbb{L}^{2}_{\nu,\mathbb{K}}$ is given by%
\begin{equation}
\left\{
\begin{array}
[c]{ll}%
dp^{1}(t) & =-[\nabla_{m}H(t)+\mathbb{E}(\nabla_{\overline{m}}H^{t}%
|\mathcal{F}_{t})]dt+q^{1}(t)dB(t)+\int_{%
\mathbb{R}
_{0}}r^{1}(t,\zeta)\widetilde{N}(dt,d\zeta);t\in\lbrack0,T],\\
p^{1}(t) & =\nabla_{m}h(X(T),M(T));t\geq T,\\
q^{1}(t) & =0;t>T,\\
r^{1}(t,\cdot) & =0;t>T,
\end{array}
\right.  \label{p1}%
\end{equation}
where $\nabla_{\bar{m}}H^{t}$ is defined in the similar way as $\nabla
_{\bar{x}}H^{t}$ above, i.e. by the property that 
\begin{equation}
{\textstyle\int_{0}^{T}}
\Gamma_{\bar{m}}(t)M(t)dt={\textstyle\int_{0}^{T}}
<\nabla_{\overline{m}}H(t),M_{t}>dt, \label{eq4.6}
\end{equation}
for all $M\in L_{0}^{2}$.


\end{itemize}

Advanced backward stochastic differential equations (absde) have been studied
by Peng and Yang \cite{PY} in the Brownian setting and for the jump case, we
refer to Øksendal et al \cite{osz}, Øksendal and Sulem \cite{ospre}. It was
also extended to the context of enlargement progressive of filtration by
Jeanblanc et al in \cite{JLA}. \newline When Agram and Røse \cite{AR} used the
maximum principle to study optimal control of mean-field delayed sfde $\left(
\ref{mfd}\right)  ,$ they obtained a mean-field absfde.

The question of existence and uniqueness of the solutions of the equations
above will not be studied here.

\subsection{A sufficient maximum principle}

We are now able to derive the sufficient version of the maximum principle.

\begin{theorem}
[Sufficient maximum principle]Let $\widehat{u}\in\mathcal{U}^{ad}$ with
corresponding solutions $\widehat{X}\in\mathcal{S}^{2}$, $(\widehat{p}%
^{0},\widehat{q}^{0},\widehat{r}^{0})\in\mathcal{S}^{2}\times\mathbb{L}%
^{2}\times\mathbb{L}_{\nu}^{2}$ and $(\widehat{p}^{1},\widehat{q}^{1}%
,\widehat{r}^{1})\in\mathcal{S}_{\mathbb{K}}^{2}\times\mathbb{L}_{\mathbb{K}%
}^{2}\times\mathbb{L}_{\nu,\mathbb{K}}^{2}$ of the forward and backward
stochastic differential equations $\left(  \ref{sfde}\right)  $, $\left(
\ref{p0}\right)  $ and $\left(  \ref{p1}\right)  $ respectively. For arbitrary
$u\in\mathcal{U}$, put
\begin{align}
&  H(t):=H(t,\widehat{X}(t),\widehat{X}_{t},\widehat{M}(t),\widehat{M}%
_{t},u(t),u_{t},\widehat{p}^{0}(t),\widehat{q}^{0}(t),\widehat{r}^{0}%
(t,\cdot),\widehat{p}^{1}(t)),\\
&  \widehat{H}(t):=H(t,\widehat{X}(t),\widehat{X}_{t},\widehat{M}%
(t),\widehat{M}_{t},\widehat{u}(t),\widehat{u}_{t},\widehat{p}^{0}%
(t),\widehat{q}^{0}(t),\widehat{r}^{0}(t,\cdot),\widehat{p}^{1}(t)).
\end{align}

Suppose that

\begin{itemize}
\item (Concavity) The functions
\[%
\begin{array}
[c]{lll}%
(x,\overline{x},m,\overline{m},u,\overline{u}) & \mapsto & H(t,x,\overline
{x},m,\overline{m},u,\overline{u},\widehat{p}^{0},\widehat{q}^{0},\widehat
{r}^{0}(\cdot),\widehat{p}^{1}),\\
(x,m) & \mapsto & h(x,m),
\end{array}
\]

are concave $\mathbb{P}$-a.s. for each $t\in\left[  0,T\right]  $.

\item (Maximum condition)%
\begin{equation}%
\begin{array}
[c]{lll}%
\mathbb{E[}\widehat{H}(t)|\mathcal{G}_{t}] & = & \underset{u\in\mathcal{U}%
^{ad}}{\sup}\mathbb{E}\left[  H(t)|\mathcal{G}_{t}\right]  ,
\end{array}
\label{maxQ}%
\end{equation}

\end{itemize}

$\mathbb{P}$-a.s. for each $t\in\left[  0,T\right]  $.

Then $\widehat{u}$ is an optimal control for the problem $\left(
\ref{perf}\right)  $.
\end{theorem}

\noindent{Proof.} \quad By considering a sequence of stopping times converging upwards to $T$, we see that we may assume that all the $dB$- and $\tilde{N}$- integrals in the following are martingales and hence have expectation 0. We refer to the proof of Lemma 3.1 in \cite{OS} for details. \\

We want to prove that $J(u)\leq J(\widehat{u})$ for
all $u\in\mathcal{U}^{ad}$. Application of definition $\left(  \ref{perf}%
\right)  $ gives for fixed $u\in\mathcal{U}^{ad}$ that
\begin{equation}%
\begin{array}
[c]{lll}%
J(u)-J(\widehat{u}) & = & I_{1}+I_{2},
\end{array}
\label{J}%
\end{equation}
where
\begin{align*}
&
\begin{array}
[c]{lll}%
I_{1} & = & \mathbb{E[}%
{\textstyle\int_{0}^{T}}
\{\ell(t)-\widehat{\ell}(t)\}dt],
\end{array}
\\
&
\begin{array}
[c]{lll}%
I_{2} & = & \mathbb{E[}h(X(T),M(T))-h(\widehat{X}(T),\widehat{M}(T))],
\end{array}
\end{align*}
with
\begin{align}
\ell(t) &  :=\ell(t,\widehat{X}(t),\widehat{X}_{t},\widehat{M}(t),\widehat
{M}_{t},u(t),u_{t}),\\
\widehat{\ell}(t) &  :=\ell(t,\widehat{X}(t),\widehat{X}_{t},\widehat
{M}(t),\widehat{M}_{t},\widehat{u}(t),\widehat{u}_{t}).
\end{align}
and similarly with $b(t),\widehat{b}(t)$ etc. later.\newline Applying the
definition of the Hamiltonian $\left(  \ref{haml}\right)  $, we get%
\begin{equation}%
\begin{array}
[c]{lll}%
I_{1} & = & \mathbb{E[}%
{\textstyle\int_{0}^{T}}
\{H(t)-\widehat{H}(t)-\widehat{p}^{0}(t)\widetilde{b}(t)-\widehat{q}%
^{0}(t)\widetilde{\sigma}(t)\\
&  & -{%
{\textstyle\int_{\mathbb{R}_{0}}}
}\widehat{r}^{0}(t,\zeta)\tilde{\gamma}(t,\zeta)\nu(d\zeta)-<\widehat{p}%
^{1}(t),\widetilde{M}^{\prime}(t)>\}dt],
\end{array}
\label{I1}%
\end{equation}
where $\widetilde{b}(t)=b(t)-\widehat{b}(t)$ etc., and
\[%
\begin{array}
[c]{lll}%
\widetilde{M}^{\prime}(t) & =\frac{d\widetilde{M}(t)}{dt} & =\frac{d}%
{dt}(M(t)-\widehat{M}(t)).
\end{array}
\]
Using concavity of $h$ and the definition of the terminal values of the absfde
$\left(  \ref{p0}\right)  $ and $\left(  \ref{p1}\right)  ,$ we get
\begin{equation}%
\begin{array}
[c]{lll}%
I_{2} & \leq & \mathbb{E}[\tfrac{\partial\widehat{h}}{\partial x}%
(T)\widetilde{X}(T)+\nabla_{m}\widehat{h}(T)\widetilde{M}(T)]\\
& = & \mathbb{E}[\widehat{p}^{0}(T)\widetilde{X}(T)+<\widehat{p}^{1}%
(T),\widetilde{M}(T)>].
\end{array}
\label{I2}%
\end{equation}
Applying the Itô formula to $\widehat{p}^{0}\widetilde{X}$ and $\widehat
{p}^{1}\widetilde{M}$, we have%
\begin{equation}%
\begin{array}
[c]{lll}%
\mathbb{E}[\widehat{p}^{0}(T)\widetilde{X}(T)] & = & \mathbb{E[}%
{\textstyle\int_{0}^{T}}
\widehat{p}^{0}(t)d\widetilde{X}(t)+%
{\textstyle\int_{0}^{T}}
\widetilde{X}(t)d\widehat{p}^{0}(t)+%
{\textstyle\int_{0}^{T}}
\widehat{q}^{0}(t)\widetilde{\sigma}(t)dt\\
&  & +%
{\textstyle\int_{0}^{T}}
{\textstyle\int_{\mathbb{R}_{0}}}
\widehat{r}^{0}(t,\zeta)\widetilde{\gamma}(t,\zeta)\nu(d\zeta)dt]\\
& = & \mathbb{E[}%
{\textstyle\int_{0}^{T}}
\widehat{p}^{0}(t)\widetilde{b}(t)dt-%
{\textstyle\int_{0}^{T}}
\tfrac{\partial\widehat{H}}{\partial x}(t)\widetilde{X}(t)dt-\int_{0}%
^{T}\mathbb{E}(\nabla_{\overline{x}}\widehat{H^{t}}|\mathcal{F}_{t}%
)\widetilde{X}(t)dt\\
&  & +%
{\textstyle\int_{0}^{T}}
\widehat{q}^{0}(t)\widetilde{\sigma}(t)dt+%
{\textstyle\int_{0}^{T}}
{\textstyle\int_{\mathbb{R}_{0}}}
\widehat{r}^{0}(t,\zeta)\widetilde{\gamma}(t,\zeta)\nu(d\zeta)dt],
\end{array}
\label{.10}%
\end{equation}
and%

\begin{align}
&\mathbb{E}[<\widehat{p}^{1}(T),\widetilde{M}(t)>]
\nonumber\\ 
& =  \mathbb{E[}%
{\textstyle\int_{0}^{T}}
<\widehat{p}^{1}(t),d\widetilde{M}(t)>dt+%
{\textstyle\int_{0}^{T}}
<\widetilde{M}(t),d\widehat{p}^{1}(t)>dt]\\
& =  \mathbb{E[}%
{\textstyle\int_{0}^{T}}
<\widehat{p}^{1}(t),\widetilde{M}^{\prime}(t)>dt -
{\textstyle\int_{0}^{T}}
<\nabla_{m}\widehat{H}(t),\widetilde{M}(t)>dt
&-{\textstyle\int_{0}^{T}}
\mathbb{E}(\nabla_{\overline{m}}\widehat{H^{t}}|\mathcal{F}_{t})\widetilde
{M}(t)dt],
\end{align}
\label{.11}
where we have used that the $dB(t)$ and $\widetilde{N}(dt,d\zeta)$ integrals
have mean zero. On substituting
$\left(  \ref{I1}\right)  ,\left(  \ref{.10}\right)  $ and $\left(
\ref{.11}\right)  $ into $\left(  \ref{J}\right)  $, we obtain%
\[%
\begin{array}
[c]{ll}%
J(u)-J(\widehat{u}) & \leq\mathbb{E[}%
{\textstyle\int_{0}^{T}}
\{H(t)-\widehat{H}(t)-%
{\textstyle\int_{0}^{T}}
\tfrac{\partial\widehat{H}}{\partial x}(t)\widetilde{X}(t)dt-%
{\textstyle\int_{0}^{T}}
\nabla_{\overline{x}}\widehat{H^{t}}\widetilde{X}(t)\}dt\\
& -%
{\textstyle\int_{0}^{T}}
<\nabla_{m}\widehat{H}(t),\widetilde{M}(t)>dt-%
{\textstyle\int_{0}^{T}}
\nabla_{\overline{m}}\widehat{H^{t}}\widetilde{M}(t)dt].
\end{array}
\]
Since $\widetilde{X}(t)=0$ for all $t\in\lbrack-\delta,0]$ and for all $t>T$
we see that $\widetilde{X}\in L_{0}^{2}$ and therefore by \eqref{eq4.6}, we
have
\begin{equation}%
{\textstyle\int_{0}^{T}}
\nabla_{\bar{x}}\widehat{H}^{t}\widetilde{X}(t)dt=%
{\textstyle\int_{0}^{T}}
<\nabla_{\bar{x}}\widehat{H}(t),\widetilde{X}_{t}>dt.\label{estpath}%
\end{equation}
Similar considerations give
\begin{equation}%
\begin{array}
[c]{ll}%
{\textstyle\int_{0}^{T}}
\nabla_{\overline{m}}\widehat{H}^{t}\widetilde{M}(t)dt & =%
{\textstyle\int_{0}^{T}}
<\nabla_{\overline{m}}\widehat{H}(t),\widetilde{M}_{t}>dt.
\end{array}
\label{stpath2}%
\end{equation}
By the assumption that $H$ is concave and that the process $u$ is
$\mathcal{G}_{t}$-adapted, we therefore get%
\[%
\begin{array}
[c]{ll}%
J(u)-J(\widehat{u}) & \leq\mathbb{E[}%
{\textstyle\int_{0}^{T}}
\{\tfrac{\partial\widehat{H}}{\partial u}(t)\widetilde{u}(t)+<\nabla
_{\overline{u}}\widehat{H}(t),\widetilde{u}_{t}>\}dt]\\
& =\mathbb{E[}%
{\textstyle\int_{0}^{T}}
\mathbb{E(}\tfrac{\partial\widehat{H}}{\partial u}(t)\widetilde{u}%
(t)+<\nabla_{\overline{u}}\widehat{H}(t),\widetilde{u}_{t}>|\mathcal{G}%
_{t})dt]\\
& =\mathbb{E[}%
{\textstyle\int_{0}^{T}}
\{\mathbb{E[}\tfrac{\partial\widehat{H}}{\partial u}(t)|\mathcal{G}%
_{t}]\widetilde{u}(t)+<\mathbb{E[}\nabla_{\overline{u}}\widehat{H}%
(t)|\mathcal{G}_{t}],\widetilde{u}_{t}>\}dt]\leq0.
\end{array}
\]
For the last inequality to hold, we use that $\mathbb{E[}\widehat
{H}(t)|\mathcal{G}_{t}]$ has a maximum at $\widehat{u}(t)$.\hfill$\square$

\subsection{A necessary maximum principle}

We now proceed to study the necessary maximum principle. Let us then impose
the following set of assumptions.

\begin{description}
\item[i)] On the coefficient functionals:
\end{description}

\begin{itemize}
\item The functions $b,\sigma$ and $\gamma$ admit bounded partial derivatives
w.r.t. $x,\overline{x},m,\overline{m},u,\overline{u}$.
\end{itemize}

\begin{description}
\item[ii)] On the performance functional:
\end{description}

\begin{itemize}
\item The function $\ell$ and the terminal value $h$ admit bounded partial
derivatives w.r.t. $x,\overline{x},m,\overline{m},u,\overline{u}$ and w.r.t.
$x,m$ respectively.
\end{itemize}

\begin{description}
\item[ii)] On the set of admissible processes:
\end{description}

\begin{itemize}
\item Whenever $u\in\mathcal{U}^{ad}$ and $\pi\in\mathcal{U}^{ad}$ is bounded,
there exists $\epsilon>0$ such that
\[
u+\lambda\pi\in\mathcal{U}^{ad}\text{, for each }\lambda\in\left[
-\epsilon,\epsilon\right]  .
\]

\item For each $t_{0}\in\left[  0,T\right]  $ and all bounded $\mathcal{G}%
_{t_{0}}$-measurable random variables $\alpha,$ the process
\[%
\begin{array}
[c]{lll}%
\pi\left(  t\right)  & = & \alpha\mathbf{1}_{\left(  t_{0},T\right]  }(t),
\end{array}
\]
belongs to $\mathcal{U}^{ad}$.\medskip
\end{itemize}

In general, if $K^{u}(t)$ is a process depending on $u$, we define the
operator $D$ on $K$ by
\begin{equation}
DK^{u}(t):=D^{\pi}K^{u}(t)=\frac{d}{d\lambda}K^{u+\lambda\pi}(t)|_{\lambda=0},
\end{equation}
whenever the derivative exists.\newline Define the \emph{derivative process}
$Z(t)$ by
\[
Z(t):=DX(t):=\tfrac{d}{d\lambda}X^{u+\lambda\pi}|_{\lambda=0}.
\]
Using matrix notation, note that $Z(t)$ satisfies the equation
\begin{equation}
\left\{
\begin{array}
[c]{ll}%
dZ\left(  t\right)   & =(\nabla b\left(  t\right)  )^{T}\text{ }(Z\left(
t\right)  ,Z_{t},DM(t),DM_{t},\pi\left(  t\right)  ,\pi_{t})dt\\
& +(\nabla\sigma\left(  t\right)  )^{T}\text{ }(Z\left(  t\right)
,Z_{t},DM(t),DM_{t},\pi\left(  t\right)  ,\pi_{t})B(t)\\
& +%
{\textstyle\int_{\mathbb{R}_{0}}}
(\nabla\gamma\left(  t,\zeta\right)  )^{T}\text{ }(Z\left(  t\right)
,Z_{t},DM(t),DM_{t},\pi\left(  t\right)  ,\pi_{t},\zeta)\widetilde
{N}(dt,d\zeta);\text{ }t\in\left[  0,T\right]  ,\\
Z\left(  t\right)   & =0;\text{ }t\in\left[  -\delta,0\right]  ,
\end{array}
\right.  \label{dervz}%
\end{equation}
where $(\nabla b)^{T}=(\tfrac{\partial b}{\partial x},\nabla_{\overline{x}%
}b,\nabla_{m}b,\nabla_{\overline{m}}b,\tfrac{\partial b}{\partial u}%
,\nabla_{\overline{u}}b)^{T}$, $(\cdot)^{T}$ denotes matrix transposed and we
mean by $\nabla_{\overline{x}}b\left(t\right)  Z_{t}$, (respectively $\nabla_{\overline
{m}}b\left(  t\right)  DM_{t})$ the action of the operator $\nabla
_{\overline{x}}b\left(  t\right)  $ $(\nabla_{\overline{m}}b\left(  t\right)
)$ on the segment $Z_{t}=\{Z(t+s)\}_{s\in\lbrack-\delta,0]}$ $(DM_{t}%
=\{DM(t+s)\}_{s\in\lbrack-\delta,0]})$ i.e., $<\nabla_{\overline{x}}b\left(
t\right)  ,Z_{t}>$ $(<\nabla_{\overline{m}}b\left(  t\right)  ,DM_{t}>)$ and
similar considerations for $\sigma$ and $\gamma$.\medskip\ 
\begin{theorem}
[Necessary maximum principle]Let $\widehat{u}\in\mathcal{U}^{ad}$ with
corresponding solutions $\widehat{X}\in\mathcal{S}^{2}$ and $(\widehat{p}%
^{0},\widehat{q}^{0},\widehat{r}^{0})\in\mathcal{S}^{2}\times\mathbb{L}%
^{2}\times\mathbb{L}^{2}_{\nu}$ and $(\widehat{p}^{1},\widehat{q}^{1}%
,\widehat{r}^{1})\in\mathcal{S}_{\mathbb{K}}^{2}\times\mathbb{L}%
^{2}_{\mathbb{K}}\times\mathbb{L}^{2}_{\nu,\mathbb{K}}$ of the forward and
backward stochastic differential equations $\left(  \ref{sfde}\right)  $ and
$\left(  \ref{p0}\right)  -\left(  \ref{p1}\right)  $ respectively, with the
corresponding derivative process $\widehat{Z}\in\mathcal{S}^{2}$ given by
$\left(  \ref{dervz}\right)  $. Then the following, (i) and (ii), are
equivalent:\medskip

\begin{description}
\item[(i)] For all bounded $\pi\in\mathcal{U}^{ad}$ 
\[%
\begin{array}
[c]{ll}%
\tfrac{d}{d\lambda}J(\hat{u}+\lambda\pi)|_{\lambda=0} & =0.
\end{array}
\]

\item[(ii)]
\[%
\begin{array}
[c]{lll}%
\mathbb{E[(}\tfrac{\partial H}{\partial u}(t)+\nabla_{\overline{u}}%
H_{t})|\mathcal{G}_{t}]_{u=\hat{u}} & = & 0\text{ for all }t\in\lbrack0,T).
\end{array}
\newline%
\]

\end{description}
\end{theorem}

\noindent{Proof.} \quad Before starting the proof, let us first clarify some
notation: Note that
\[%
\begin{array}
[c]{ll}%
\nabla_{m}<p_{1}^{1}(t),\frac{d}{dt}m> & =<p_{1}^{1}(t),\frac{d}{dt}(\cdot)>,
\end{array}
\]
and hence
\[%
\begin{array}
[c]{llll}%
<\nabla_{m}<p_{1}^{1}(t),\frac{d}{dt}m>,DM(t)> & =<p_{1}^{1}(t),\frac{d}%
{dt}DM(t)> & =<p_{1}^{1}(t),DM^{\prime}(t)> & =p_{1}^{1}(t)DM^{\prime}(t).
\end{array}
\]
Also, note that%

\begin{equation}
dDM(t)=DM^{\prime}(t)dt.
\end{equation}
\quad By considering a sequence of stopping times converging upwards to $T$, we see that we may assume that all the $dB$- and $\tilde{N}$- integrals in the following are martingales and hence have expectation 0. We refer to the proof of Lemma 3.1 in \cite{OS} for details. \\

Assume that (i) holds. Then
\begin{align}
0 & =\tfrac{d}{d\lambda}J(u+\lambda\pi)|_{\lambda=0}=\mathbb{E[}
{\textstyle\int_{0}^{T}}
\{(\nabla\ell\left(  t\right)  )^{T}\text{ }(Z\left(  t\right)  ,Z_{t}
,DM(t),DM_{t},\pi\left(  t\right)  ,\pi_{t})\}dt\\
& +\tfrac{\partial h}{\partial x}\left(  T\right)  Z\left(  T\right)
+\nabla_{m}h\left(  T\right)  DM(T)].
\end{align}

Hence, by the definition of $H$ $\left(  \ref{haml}\right)  $ and the terminal
values of the absfde $p^{0}(T)$ and $p^{1}(T),$\ we have
\[%
\begin{array}
[c]{ll}%
0 & =\tfrac{d}{d\lambda}J(u+\lambda\pi)|_{\lambda=0}\\
& =\mathbb{E[}%
{\textstyle\int_{0}^{T}}
\{(\nabla H\left(  t\right)  )^{T}\text{ }(Z\left(  t\right)  ,Z_{t}%
,DM(t),DM_{t},\pi\left(  t\right)  ,\pi_{t})\\
& -p^{0}(t)(\nabla b\left(  t\right)  )^{T}\text{ }(Z\left(  t\right)
,Z_{t},DM(t),DM_{t},\pi\left(  t\right)  ,\pi_{t})\\
& -q^{0}(t)(\nabla\sigma\left(  t\right)  )^{T}\text{ }(Z\left(  t\right)
,Z_{t},DM(t),DM_{t},\pi\left(  t\right)  ,\pi_{t})\\
& -%
{\textstyle\int_{\mathbb{R}_{0}}}
r^{0}(t,\zeta)(\nabla\gamma\left(  t,\zeta\right)  )^{T}\text{ }(Z\left(
t\right)  ,Z_{t},DM(t),DM_{t},\pi\left(  t\right)  ,\pi_{t})\nu(d\zeta)\}dt]\\
& -%
{\textstyle\int_{0}^{T}}
p^{1}(t)DM^{\prime}\left(  t\right)  dt+p^{0}(T)Z(T)+p^{1}(T)DM(T)].
\end{array}
\]
Applying Itô formula to both $p^{0}Z$ and $p^{1}DM$, we get%
\[%
\begin{array}
[c]{ll}%
\mathbb{E}[p^{0}(T)Z(T)] & =\mathbb{E[}%
{\textstyle\int_{0}^{T}}
p^{0}(t)dZ(t)+%
{\textstyle\int_{0}^{T}}
Z(t)dp^{0}(t)+[p^{0},Z]_{T}]\\
& =\mathbb{E[}%
{\textstyle\int_{0}^{T}}
p^{0}(t)(\nabla b\left(  t\right)  )^{T}\text{ }(Z\left(  t\right)
,Z_{t},DM(t),DM_{t},\pi\left(  t\right)  ,\pi_{t})dt\\
& -%
{\textstyle\int_{0}^{T}}
\{\tfrac{\partial H}{\partial x}(t)+\nabla_{\overline{x}}H^{t}\}Z(t)dt\\
& +%
{\textstyle\int_{0}^{T}}
q^{0}(t)(\nabla\sigma\left(  t\right)  )^{T}\text{ }(Z\left(  t\right)
,Z_{t},DM(t),DM_{t},\pi\left(  t\right)  ,\pi_{t})dt\\
& +%
{\textstyle\int_{0}^{T}}
{\textstyle\int_{\mathbb{R}_{0}}}
r^{0}(t,\zeta)(\nabla\gamma\left(  t,\zeta\right)  )^{T}\text{ }(Z\left(
t\right)  ,Z_{t},DM(t),DM_{t},\pi\left(  t\right)  ,\pi_{t})\nu(d\zeta)dt],
\end{array}
\]
and
\begin{equation}%
\begin{array}
[c]{ll}%
\mathbb{E[}<p^{1}(T),DM(T)>] & =\mathbb{E[}%
{\textstyle\int_{0}^{T}}
<p^{1}(t),DM^{\prime}(t)>dt+%
{\textstyle\int_{0}^{T}}
<dp^{1}(t),DM(t)>]\\
& =\mathbb{E[}%
{\textstyle\int_{0}^{T}}
<p^{1}(t),DM^{\prime}(t)>dt-%
{\textstyle\int_{0}^{T}}
< \{\nabla_{m}H(t)+\nabla_{\overline{m}}H^{t}\},DM(t)>dt].
\end{array}
\nonumber
\end{equation}
Proceeding as in $\left(  \ref{estpath}\right)  -\left(  \ref{stpath2}\right)
,$ we obtain%
\[%
\begin{array}
[c]{ll}%
{\textstyle\int_{0}^{T}}
\nabla_{\overline{x}}H^{t}Z(t)dt & =%
{\textstyle\int_{0}^{T}}
<\nabla_{\overline{x}}H(t),Z_{t}>dt,\\%
{\textstyle\int_{0}^{T}}
\nabla_{\overline{m}}H^{t}DM(t)dt & =%
{\textstyle\int_{0}^{T}}
<\nabla_{\overline{m}}H(t),DM_{t}>dt.
\end{array}
\]
Combining the above, we get%

\begin{equation}%
\begin{array}
[c]{ll}%
0 & =\mathbb{E[}%
{\textstyle\int_{0}^{T}}
(\tfrac{\partial H}{\partial u}(t)\pi(t)+\left\langle \nabla_{\overline{u}%
}H(t),\pi_{t}\right\rangle )dt]\text{.}%
\end{array}
\label{h_pi}%
\end{equation}
Now choose $\pi(t)=\alpha\mathbf{1}_{\left(  t_{0},T\right]  }(t)$, where
$\alpha=\alpha(\omega)$ is bounded and $\mathcal{G}_{t_{0}}$-measurable and
$t_{0}\in\lbrack0,T)$. Then $\pi_{t}=\alpha\{\mathbf{1}_{\left(
t_{0},T\right]  }(t+s)\}_{s\in\lbrack-\delta,0]}$ and $\left(  \ref{h_pi}%
\right)  $ gives
\[%
\begin{array}
[c]{ll}%
0 & =\mathbb{E[}%
{\textstyle\int_{t_{0}}^{T}}
\tfrac{\partial H}{\partial u}(t)\alpha dt+%
{\textstyle\int_{t_{0}}^{T}}
\left\langle \nabla_{\overline{u}}H(t),\alpha\right\rangle \{\mathbf{1}%
_{\left(  t_{0},T\right]  }(t+s)\}_{s\in\lbrack-\delta,0]}dt]\text{.}%
\end{array}
\]
Differentiating with respect to $t_{0},$ we obtain%
\[%
\begin{array}
[c]{ll}%
\mathbb{E[}(\tfrac{\partial H}{\partial u}(t_{0})+\nabla_{\overline{u}%
}H_{t_{0}})\alpha] & =0\text{,}%
\end{array}
\]
Since this holds for all such $\alpha$, we conclude that%
\[%
\begin{array}
[c]{ll}%
\mathbb{E[}(\tfrac{\partial H}{\partial u}(t_{0})+\nabla_{\overline{u}%
}H_{t_{0}})|\mathcal{G}_{t_{0}}] & =0\text{, which is (ii).}%
\end{array}
\]
This argument can be reversed, to prove that (ii)$\Longrightarrow$(i). We omit
the details.\hfill$\square$ \newline

\section{Applications}

We illustrate our results by studying some examples.

\subsection{Mean-variance portfolio with memory}

We apply the results obtained in the previous sections to solve the memory
mean-variance problem by proceeding as it has been done in Framstad et al
\cite{FOS}, Anderson and Djehiche \cite{AD} and Røse \cite{R}.$\newline%
$Consider the state equation $X^{\pi}(t)=X(t)$ on the form%
\begin{equation}
\left\{
\begin{array}
[c]{ll}%
dX(t) & =X(t-\delta)\pi(t)[b_{0}(t)dt+\sigma_{0}(t)dB(t)+%
{\textstyle\int_{\mathbb{R}_{0}}}
\gamma_{0}\left(  t,\zeta\right)  \widetilde{N}(dt,d\zeta)];t\in\left[
0,T\right]  ,\\
X(t) & =\xi(t);t\in\left[  -\delta,0\right]  ,
\end{array}
\right.  \label{w}%
\end{equation}

\noindent for some bounded deterministic function $\xi(t);t\in\lbrack
-\delta,0]$. We assume that the admissible processes are càdlàg processes in
$L^{2}(\Omega,[0,T])$, that are adapted to the filtration $\mathcal{F}_{t}$
and such that a unique solution exists. The coefficients $b_{0},\sigma_{0}$
and $\gamma_{0}>-1$ are supposed to be bounded $\mathbb{F}$-adapted processes
with
\[
|b_{0}(t)|>0\text{ and }\sigma_{0}^{2}(t)+%
{\textstyle\int_{\mathbb{R}_{0}}}
\gamma_{0}^{2}(t,\zeta)\nu(d\zeta)>0\text{ a.s. for all }t.
\]
We want to find an admissible portfolio $\pi(t)$ which maximizes
\begin{equation}%
\begin{array}
[c]{lll}%
J(\pi) & = & \mathbb{E[}-\frac{1}{2}(X(T)-a)^{2}],
\end{array}
\label{p}%
\end{equation}
over the set of admissible processes $\mathcal{U}^{ad}$ and for a given
constant $a\in%
\mathbb{R}
$.$\newline$The Hamiltonian for this problem is given by%
\begin{equation}%
\begin{array}
[c]{lll}%
H(t,\overline{x},\pi,p^{0},q^{0},r^{0}(\cdot)) & = & \pi G(\overline{x}%
)(b_{0}p^{0}+\sigma_{0}q^{0}+%
{\textstyle\int_{\mathbb{R}_{0}}}
\gamma_{0}\left(  \zeta\right)  r^{0}\left(  \zeta\right)  \nu(d\zeta)),
\end{array}
\label{h}%
\end{equation}
where
\begin{equation}
G(\bar{x})=x(\delta)\text{ when }\bar{x}=\{x(s)\}_{s\in [0,\delta]},
\end{equation}
i.e. $G$ is evaluation at $r=\delta.$
See Example 4.4 (i). Hence by Lemma 4.3 the triple $(p^{0},q^{0},r^{0}%
)\in\mathcal{S}^{2}\times\mathbb{L}^{2}\times\mathbb{L}_{\nu}^{2}$ is the
adjoint process which satisfies
\begin{equation}
\left\{
\begin{array}
[c]{ll}%
dp^{0}(t) & =-\mathbb{E[}\pi(t+\delta)(b_{0}(t+\delta)p^{0}(t+\delta
)+\sigma_{0}(t+\delta)q^{0}(t+\delta)\\
& +%
{\textstyle\int_{\mathbb{R}_{0}}}
\gamma_{0}\left(  t+\delta,\zeta\right)  r^{0}\left(  t+\delta,\zeta\right)
\nu(d\zeta))|\mathcal{F}_{t}]dt+q^{0}(t)dB(t)\\
& +%
{\textstyle\int_{\mathbb{R}_{0}}}
r^{0}\left(  t,\zeta\right)  \widetilde{N}(dt,d\zeta);t\in\left[  0,T\right]
,\\
p^{0}(t) & =-(X(T)-a);t\geq T,\\
q^{0}(t) & =r^{0}(\cdot)=0;t>T.
\end{array}
\right.  \label{ap}%
\end{equation}
Existence and uniqueness of equations of type $\left(  \ref{ap}\right)  $ have
been studied by Øksendal et al \cite{osz}.\newline

Suppose that $\widehat{\pi}$ is an optimal control. Then by the necessary
maximum principle, we get for each $t$ that%
\begin{align}
0  &  =\tfrac{\partial\widehat{H}}{\partial\pi}(t,\widehat{X_{t}},\widehat
{\pi}(t),\widehat{p}^{0}(t),\widehat{q}^{0}(t),\widehat{r}^{0}(t,\cdot
))\label{nc}\\
&  =\widehat{X}(t-\delta)(b_{0}(t)\widehat{p}^{0}(t)+\sigma_{0}(t)\widehat
{q}^{0}(t)+%
{\textstyle\int_{\mathbb{R}_{0}}}
\gamma_{0}\left(  t,\zeta\right)  \widehat{r}^{0}\left(  t,\zeta\right)
\nu(d\zeta)).\nonumber
\end{align}
So we search for a candidate $\widehat{\pi}$ satisfying%
\begin{equation}%
\begin{array}
[c]{ll}%
0 & =b_{0}(t)\widehat{p}^{0}(t)+\sigma_{0}(t)\widehat{q}^{0}(t)+%
{\textstyle\int_{\mathbb{R}_{0}}}
\gamma_{0}\left(  t,\zeta\right)  \widehat{r}^{0}\left(  t,\zeta\right)
\nu(d\zeta)\text{, for all }t.
\end{array}
\label{pi}%
\end{equation}
This gives the following adjoint equation:
\begin{equation}
\left\{
\begin{array}
[c]{ll}%
d\widehat{p}^{0}(t) & =\widehat{q}^{0}(t)dB(t)+%
{\textstyle\int_{\mathbb{R}_{0}}}
\widehat{r}^{0}\left(  t,\zeta\right)  \widetilde{N}(dt,d\zeta);t\in\left[
0,T\right]  ,\\
\widehat{p}^{0}(t) & =-(X(T)-a);t\geq T,\\
\widehat{q}^{0}(t) & =\widehat{r}^{0}(\cdot)=0;t>T.
\end{array}
\right.  \label{ap1}%
\end{equation}
We start by guessing that $\widehat{p}^{0}$ has the form{\normalsize
\begin{equation}%
\begin{array}
[c]{lll}%
\widehat{p}^{0}(t) & = & \varphi(t)\widehat{X}(t)+\psi(t)
\end{array}
\label{p^}%
\end{equation}
}for some deterministic functions $\varphi${\normalsize $,\psi\in C^{1}[0,T]$}
with{\normalsize
\begin{equation}
\varphi(T)=-1,\quad\psi(T)=a. \label{v}%
\end{equation}
}Using the Itô formula to find the integral representation of $\widehat{p}%
^{0}$ and comparing with the adjoint equation (\ref{ap1}), we find that the
following three equations need to be satisfied:
\begin{equation}%
\begin{array}
[c]{lll}%
0 & = & \varphi^{\prime}(t)\widehat{X}(t)+\psi^{\prime}(t)+\varphi
(t)\widehat{X}(t-\delta)\widehat{\pi}(t)b_{0}(t),
\end{array}
\label{d}%
\end{equation}%
\begin{equation}%
\begin{array}
[c]{lll}%
\widehat{q}^{0}(t) & = & \varphi(t)\widehat{X}(t-\delta)\widehat{\pi}%
(t)\sigma_{0}(t),
\end{array}
\label{di}%
\end{equation}%
\begin{equation}%
\begin{array}
[c]{lll}%
\widehat{r}^{0}(t,\zeta) & = & \varphi(t)\widehat{X}(t-\delta)\widehat{\pi
}(t)\gamma_{0}(t,\zeta).
\end{array}
\label{j}%
\end{equation}
Assuming that $\widehat{X}${\normalsize $(t)\neq0$ }$\mathbb{P}%
{\normalsize \times dt}$-a.e. and $\varphi(t)\neq0$ for each $t$, we find from
equation $(\ref{d})$ that $\widehat{\pi}$ needs to satisfy{\normalsize
\[%
\begin{array}
[c]{lll}%
\widehat{\pi}(t) & = & -\frac{\varphi^{\prime}(t)\widehat{X}(t)+\psi^{\prime
}(t)}{\varphi(t)\widehat{X}(t-\delta)b_{0}(t)}\text{.}%
\end{array}
\]
}Now inserting the expressions for the adjoint processes $(\ref{d})$,
$(\ref{di})$ and $(\ref{j})$ into $(\ref{pi})$, the following equation need to
be satisfied:{\small
\[
0=b_{0}(t)[\varphi(t)\widehat{X}(t)+\psi(t)]+\varphi(t)\widehat{X}%
(t-\delta)\widehat{\pi}(t)\big(\sigma_{0}^{2}(t)+%
{\textstyle\int_{\mathbb{R}_{0}}}
\gamma_{0}^{2}(t,\zeta)\nu(d\zeta)\big).
\]
}This means that the control $\widehat{\pi}$ also needs to
satisfy{\normalsize
\begin{equation}%
\begin{array}
[c]{lll}%
\widehat{\pi}(t) & = & -\tfrac{b_{0}(t)[\varphi(t)\widehat{X}(t)+\psi
(t)]}{[\sigma_{0}^{2}(t)+%
{\textstyle\int_{\mathbb{R}_{0}}}
\gamma_{0}^{2}(t,\zeta)\nu(d\zeta)]\varphi(t)\widehat{X}(t-\delta)}.
\end{array}
\label{pih}%
\end{equation}
}By comparing the two expressions for $\widehat{\pi}$, we find that{\small
\begin{equation}%
\begin{array}
[c]{l}%
b_{0}^{2}(t)[\varphi(t)\widehat{X}(t)+\psi(t)]\\
=(\sigma_{0}^{2}(t)+%
{\textstyle\int_{\mathbb{R}_{0}}}
\gamma_{0}^{2}(t,\zeta)\nu(d\zeta))[\varphi^{\prime}(t)\widehat{X}%
(t)+\psi^{\prime}(t)].
\end{array}
\label{co}%
\end{equation}
}Now define{\normalsize
\begin{equation}%
\begin{array}
[c]{lll}%
\Lambda(t) & := & \tfrac{b_{0}^{2}(t)}{\sigma_{0}^{2}(t)+\int_{\mathbb{R}_{0}%
}\gamma_{0}^{2}(t,\zeta)\nu(d\zeta)}.
\end{array}
\label{lam}%
\end{equation}
}Then from equation $(\ref{co})$, we need to have{\normalsize
\[%
\begin{array}
[c]{lll}%
\varphi^{\prime}(t)-\Lambda(t)\varphi(t) & = & 0,\\
\psi^{\prime}(t)-\Lambda(t)\psi(t) & = & 0.
\end{array}
\]
}Together with the terminal values $(\ref{v})$, these equations have the
solution{\normalsize
\[%
\begin{array}
[c]{lll}%
\varphi(t) & = & -\exp(-%
{\textstyle\int_{t}^{T}}
\Lambda(s)ds),\\
\psi(t) & = & a\exp(-%
{\textstyle\int_{t}^{T}}
\Lambda(s)ds).
\end{array}
\]
}Then from equation $(\ref{pih})$ we can compute

{\normalsize
\[
\widehat{\pi}(t)=\tfrac{b_{0}(t)\left(  \widehat{X}(t)-\tfrac{\psi(t)}%
{\varphi(t)}\right)  }{\sigma_{0}^{2}(t)+%
{\textstyle\int_{\mathbb{R}_{0}}}
\gamma_{0}^{2}(t,\zeta)\nu(d\zeta))\widehat{X}(t-\delta)}=\Lambda
(t)\tfrac{\left(  \widehat{X}(t)-\tfrac{\psi(t)}{\varphi(t)}\right)  }%
{b_{0}(t)\widehat{X}(t-\delta)}=\tfrac{\Lambda(t)}{b_{0}(t)\widehat
{X}(t-\delta)}(\widehat{X}(t)-a).
\]
}Now, with our choice of $\widehat{\pi}$, the corresponding state equation is
the solution of%
\begin{equation}
\left\{
\begin{array}
[c]{ll}%
d\widehat{X}(t) & =\tfrac{\Lambda(t)}{b_{0}(t)}(\widehat{X}(t)-a)[b_{0}%
(t)dt+\sigma_{0}(t)dB(t)+%
{\textstyle\int_{\mathbb{R}_{0}}}
\gamma_{0}\left(  t,\zeta\right)  \widetilde{N}(dt,d\zeta)];t\in\left[
0,T\right]  ,\\
\widehat{X}(t) & =x_{0}(t);t\in\left[  -\delta,0\right]  .
\end{array}
\right.  \label{wp}%
\end{equation}
Put $Y(t)=\widehat{X}(t)-a$, then%
\begin{equation}%
\begin{array}
[c]{ll}%
dY(t) & =Y(t)[\Lambda(t)b_{0}(t)dt+\tfrac{\Lambda(t)}{b_{0}(t)}\sigma
_{0}(t)dB(t)+\int_{\mathbb{R}_{0}}\tfrac{\Lambda(t)}{b_{0}(t)}\gamma
_{0}\left(  t,\zeta\right)  \widetilde{N}(dt,d\zeta)].
\end{array}
\label{lw}%
\end{equation}
The linear equation $(\ref{lw})$ has the following explicit solution%
\[%
\begin{array}
[c]{ll}%
Y(t) & =Y(0)\exp[\int_{0}^{t}\Lambda(s)b_{0}(s)ds+\int_{0}^{t}\tfrac
{\Lambda(s)}{b_{0}(s)}\sigma_{0}(s)dB(s)+\int_{0}^{t}\int_{\mathbb{R}_{0}%
}\tfrac{\Lambda(s)}{b_{0}(s)}\gamma_{0}\left(  s,\zeta\right)  \widetilde
{N}(ds,d\zeta)].
\end{array}
\]
So if $Y(0)>0$ then $Y(t)>0$ for all $t$.\newline We have proved the following:

\begin{theorem}
[Optimal mean-variance portfolio]Suppose that $\xi(t)>a$ for all $t\in\left[
-\delta,0\right]  $. Then $\widehat{X}(t-\delta)>0$ for all $t\geq0$ and the
solution $\widehat{\pi}\in\mathcal{U}^{ad}$ of the mean-variance portfolio
problem \eqref{p} is given in feedback form as{\normalsize
\[%
\begin{array}
[c]{lll}%
\widehat{\pi}(t) & = & \tfrac{\Lambda(t)}{b_{0}(t)\widehat{X}(t-\delta
)}(\widehat{X}(t)-a),
\end{array}
\]
}where $\widehat{X}(t)$ and $\Lambda(t)$ are given by equations $(\ref{wp})$
and $(\ref{lam})$ respectively.
\end{theorem}

\subsection{A linear-quadratic (LQ) problem with memory}

We now consider a linear-quadratic control problem for a controlled system
$X(t)=X^{u}(t)$ driven by a distributed delay, of the form
\begin{equation}
\left\{
\begin{array}
[c]{ll}%
dX(t) & =[%
{\textstyle\int_{0}^{\delta}}
a(s)X(t-s)ds+u(t)]dt+\alpha_{0}(t)dB(t)+\int_{\mathbb{R}_{0}}\beta_{0}\left(t,
\zeta\right)  \widetilde{N}(dt,d\zeta);t\in\left[  0,T\right]  ,\\
X(t) & =\xi(t);t\in\lbrack-\delta,0],
\end{array}
\right.  \label{f}%
\end{equation}
where $\xi(\cdot)$ and $a(\cdot)$ are given bounded deterministic functions,
$\alpha_{0}(\cdot)$ and $\beta_{0}(\cdot,\zeta)$ are given bounded predictable processes
and $u\in\mathcal{U}^{ad}$ is our control process. We want to minimize the
expected value of $X^{2}(T)$ with a minimal average use of energy, measured by
the integral $\mathbb{E}[%
{\textstyle\int_{0}^{T}}
u^{2}(t)dt]$, i.e. the performance functional is of the quadratic type
\[
J(u)=-\tfrac{1}{2}\mathbb{E}[X^{2}(T)+%
{\textstyle\int_{0}^{T}}
u^{2}(t)dt].
\]
Our goal is to find $\widehat{u}\in\mathcal{U}^{ad},$ such that%
\begin{equation}
J(\widehat{u})=\underset{u\in\mathcal{U}^{ad}}{\sup}J(u). \label{eq5.20}%
\end{equation}
The Hamiltonian in that case takes the form
\[
H(t,\overline{x},u,p^{0},q^{0},r^{0})=-\tfrac{1}{2}u^{2}+(F(\overline{x}%
)+u)p^{0}+\alpha_{0}(t)q^{0}+%
{\textstyle\int_{\mathbb{R}_{0}}}
r^{0}\left(  \zeta\right)  \beta_{0}(t,\zeta)  \nu(d\zeta),
\]
where
\begin{equation}
F(\bar{x})=%
{\textstyle\int_{0}^{\delta}}
a(s)x(s)ds\text{ when }\bar{x}=\{x(s)\}_{s\in [0,\delta]}.
\end{equation}
By Lemma 4.3 and Example 4.4 (i) we see that the adjoint absde for
$(p^{0},q^{0},r^{0})$ is the following linear absde
\begin{equation}
\left\{
\begin{array}
[c]{ll}%
dp^{0}(t) & =-\mathbb{E}[%
{\textstyle\int_{0}^{\delta}}
a(r)p^0(t+r)dr|\mathcal{F}_{t}]dt+q^{0}(t)dB(t)+%
{\textstyle\int_{\mathbb{R}_{0}}}
r^{0}\left(  t,\zeta\right)  \widetilde{N}(dt,d\zeta);\hspace{0.2cm}%
t\in\left[  0,T\right]  ,\\
p^{0}(T) & =-X(T);t\geq T.
\end{array}
\right.  \label{b}%
\end{equation}
The function $u\mapsto H(t,\widehat{X}(t-\delta),u,\widehat
{p}^{0}(t),\widehat{q}^{0}(t),r^{0}\left(  t,\zeta\right)  )$ is maximal when%
\begin{equation}
u(t)=\widehat{u}(t)=\widehat{p}^{0}(t). \label{eq5.23}%
\end{equation}
We have proved:

\begin{theorem}
The optimal control $\hat{u}$ of the LQ memory problem \eqref{eq5.20} is given
by \eqref{eq5.23}, where the quadruplet $(\hat{X}(t)=X^{\hat{u}}(t),\hat
{p}^{0}(t),\hat{q}^{0}(t),\hat{r}^{0}\left(  t,\zeta\right)  )$ solves the following coupled
system of forward-backward stochastic
differential equations with distributed delay:

\begin{itemize}
\item
\begin{align}\label{eq5.24}
&d\hat{X}(t)=\Big(  {\textstyle\int_{0}^{\delta}}
a(r)\hat{X}(t-r)dr+\hat{p}^0(t) \Big )dt+\alpha_{0}(t)dB(t)\nonumber+\int_{\mathbb{R}_{0}}\beta_{0}\left(t,
\zeta\right)  \widetilde{N}(dt,d\zeta);t\in\left[  0,T\right] ,\nonumber\\
&\hat{X}(t)  =\xi(t);t\in\lbrack-\delta,0],
\end{align}
\item
\begin{align}\label{eq5.25}
&d\widehat{p}^{0}(t) =-\Big({\textstyle\int_{0}^{\delta}}
a(r)\mathbb{E}[\widehat{p}^0(t+r)|\mathcal{F}_{t}]dr \Big) dt+\widehat{q}^{0}(t)dB(t)+{\textstyle\int_{\mathbb{R}_{0}}}
\widehat{r}^{0}\left(  t,\zeta\right)  \widetilde{N}(dt,d\zeta);\hspace
{0.2cm}t\in\left[  0,T\right]  ,\nonumber\\
&\widehat{p}^{0}(T) =-\widehat{X}(T);t\geq T.
\end{align}

\end{itemize}
\end{theorem}

\begin{remark}
We may regard this coupled system \eqref{eq5.24}-\eqref{eq5.25} as the
corresponding Riccati equation to our LQ memory problem. See e.g. Hu \& \O ksendal \cite{HO},
page 1747.
\end{remark}

\section{Acknowledgments}
We want to thank Rosestolato Mauro for helpful comments.

\end{document}